\documentclass[11pt]{article}

\title{$q$-Poisson, $q$-Dobinski, $q$-Rota and $q$-coherent states-
 a second fortieth anniversary  memoir }
\author{A.K.Kwa\'sniewski\\  
\\ High School of Mathematics and Applied Informatics\\
PL - 15-021 Bialystok , ul.Kamienna 17,  Poland
\\e-mail: kwandr@uwb.edu.pl}

\usepackage{amsmath,amsthm}

\chardef\bslash=`\\ 
\begin{document}
\maketitle \ ArXiv: math.CO/0402254 vol 1 16 Feb 2004  corrected
\begin{abstract}
The $q$- Dobinski formula may  be interpreted as the average of
powers of random variable $X_q$ with the $q$-Poisson distribution.
\end{abstract}

Forty years ago Rota G. C. \cite{1}  proved the exponential
generating function for Bell numbers $B_n$ to be of the form

\begin{equation}\label{eq1}
     \sum_{n=0}^\infty \frac {x^n}{n!}(B_n)= \exp(e^x-1)
\end{equation}
 using the linear functional  \textit{L } such that
\begin{equation}\label{eq2}
     L(X^{\underline{n}})=1  , \qquad  n\geq 0
\end{equation}\label{eq3}
Then Bell numbers (see: formula (4)  in \cite{1})  are defined by
\begin{equation}\label{eq4}
          L(X^n)=B_n  ,\qquad   n \geq 0
\end{equation}
The above formula is exactly the Dobinski formula \cite{2} if $L$
is interpreted as the average functional for the random variable
$X$  with the Poisson distribution with $L(X) = 1$. It is Blissard
calculus inspired umbral formula \cite{1} .
    Recently an interest to Stirling numbers and consequently
to Bell numbers was revived among "$q$-coherent states physicists"
\cite{3,4,5}. Namely the expectation value with respect to
coherent state $|\gamma >$ with $|\gamma| = 1$ of the $n$-th power
of the number of quanta operator is "just" the $n$-th Bell number
$B_n$ and the explicit formula for this expectation number of
quanta is "just" Dobinski formula \cite{3}. The same is with the
$q$-coherent states case \cite{3} i.e. the expectation value with
respect to $q$-coherent state $|\gamma> $ with  $|\gamma| = 1$ of
the $n$-th power of the number operator is the $n$-th  $q$-Bell
number \cite{6,3} defined as the sum  of $q$-Stirling numbers  $
\Big\{{n \atop; k}\Big\}_q $  due to Carlitz as in \cite{6,3,4,5}.
Note there then that  for  standard Gauss $q$-extension $ x_q $ of
number  $ x $ we have
\begin{equation}\label{eq5}
x_q^n=\sum_{k=0}^{n}\Big\{{n \atop k}\Big\}_q  x_q^{\underline k}
\end{equation}
Hence the expectation value with respect to $q$-coherent state $
|\gamma> $ with  $ |\gamma|=1 $ of the  $n$-th power of the number
operator  is exactly the popular $q$-Dobinski formula . It can be
given via (3) Blissard calculus inspired umbral formula form and
may be treated as definition of $ B_n(q)$
\begin{equation}\label{eq6}
         L_q(X_q^n)=B_n(q)  ,\qquad n\geq 0 .
\end{equation}
due to the fact that  linear functional  $L_q$    interpreted as
the average functional for the random variable  $X_q$  with the
$q$-Poisson distribution with  $L_q(X_q )= 1$  satisfies
\begin{equation}\label{eq7}
L_q(X_q^{\underline{n}})=1  , \qquad  n\geq 0.
\end{equation}
We arrive to this simple conclusion using Jackson derivative
difference operator in place of  $D = d/dx$  in $q$ =1 case and
the power series generating function $G(t)$  for $q$-Poisson
probability distribution:
\begin{equation}\label{eq8}
p_k=[exp_q\lambda]^{-1}\frac{\lambda^k}{k_q!},
G(t)=\sum_{n\geq0}p_k t^k ,
\end{equation}
\begin{equation}\label{eq9}
p_n=[\frac{\partial_q^n G(t)}{n_q!}]_{t=0} , [\partial_q
G(t)]_{t=1} = 1  for   \lambda = 1.
\end{equation}
There are many $q$-extensions of Stirling numbers according to
their weighted counting interpretation. For example $ w(\pi) =
q^{cross(\pi)} , w(\pi) = q^{inv(\pi)}$  from \cite{7} gives after
 being summed over  the  set of  $k$-block partitions  the
 Carlitz $q$-Stirling numbers  or $w(\pi) = q^{nin(\pi)} $
from \cite{8} gives rise to Carlitz-Gould $q$- Stirling numbers
after being summed over  the  set of  $k$-block partitions
or with  $w(\pi) = q^{i(\pi)}$ in \cite{9} - we arrive at another
combinatorial interpretation of $q$-extended Stirling numbers.
$q$-Stirling numbers much different from Carlitz $q$-ones were
introduced in the reference \cite{10} from where one infers \cite{11}
the \textit{cigl}-analog of (5) .  Let $\Pi$ denotes  the lattice of all
partitions of the set $\{0,1,..,n-1\}$.Let  $\pi\in \Pi $  be
represented by blocks   $ \pi =\{B_o ,B_1 ,...B_i ,...\}$ , where
$B_o$ is the block containing zero:  $0\in B_o$. The weight
adapted by Cigler  defines  weighted partitions` counting
according to the content of $B_o$. Namely $w(\pi)= q^{cigl(\pi)},
cigl(\pi)=\sum_{l\in {B_0}}l $,
$\sum_{\pi\in A_{n,k}}{q^{cigl(\pi)}}\equiv\Big\{{n\atop k}\Big\}_q $
therefore $ \sum_{\pi\in \Pi}{q^{cigl(\pi)}}\equiv{B_n(q)}$.
Here $A_{n,k}$ stays for subfamily of all $k$-block partitions.
With the above relations one has defined the \textit{cigl}-$q$-Stirling
 and the \textit{cigl}-$q$-Bell numbers.
The \textit{cigl}-$q$-Stirling numbers of the second kind  are
expressed in terms of $q$-binomial coefficients and $q =1$
Stirling numbers of the second kind \cite{10} . These are new
$q$-Stirling numbers. The corresponding \textit{cigl}-$q$-Bell
numbers recently have been equivalently defined via
\textit{cigl}-$q$-Dobinski formula \cite{11} $L(X_q^n)=B_n(q),
\qquad n\geq 0 , \qquad X_q^n\equiv X(X+q-1)...(X-1+q^{n-1}) $
interpreted as  the average of this specific  $n-th$
\textit{cigl}-$q$-power random variable $ X_q^{n} $ with the $q =
1$ Poisson distribution such that $ L(X)=1 $. To this end note
that in \cite{12} , \cite{13} a family of  the so called
$\psi$-Poisson processes was introduced. The corresponding choice
of the function sequence $\psi$ leads  to the $q$-Poisson process.
Accordingly the extension of Dobinski formula with its elementary
essential content and context to general case of $\psi$- umbral
instead $q$-umbral calculi case only - is automatic in view of  an
experience from \cite{12} , \cite{13} (see corresponding earlier
references there and necessary definitions). At first what you do
is to replace index $q$ by $\psi$ in formulas (3), (4),...,(8) .
Then you have got started problems with not easy combinatorial
interpretation if at all and... etc. $\psi$-Stirling numbers and
$\psi$-Bell numbers are being then defined by (4) and (3)
correspondingly with $q$ replaced by $\psi$. We get used to write
these extensions in mnemonic convenient
 upside down notation  \cite{12} , \cite{13}
\begin{equation}\label{eq10}
\psi_n\equiv n_\psi , x_{\psi}\equiv \psi(x)\equiv\psi_x ,
 n_\psi!=n_\psi(n-1)_\psi!, n>0 ,
\end{equation}
\begin{equation}\label{eq11}
x_{\psi}^{\underline{k}}=x_{\psi}(x-1)_\psi(x-2)_{\psi}...(x-k+1)_{\psi}
\end{equation}
\begin{equation}\label{eq12}
x_{\psi}(x-1)_{\psi}...(x-k+1)_{\psi}= \psi(x)
\psi(x-1)...\psi(x-k-1) .
\end{equation}
You may consult for further development and use of this notation
\cite{12} , \cite{13} and references therein.


\begin{thebibliography}{99}
\parskip 0pt

\bibitem{1}
Rota G. C. {\it The number of partitions of a set} Amer. Math.
Monthly {\bf 71}(1964) :  498-504

\bibitem{2}
G. Dobinski {\it Summierung der Reihe S .... für m = 1, 2, 3, 4,
5, ....} Grunert Archiv (Arch. Math. Phys) {\bf 61}, 333-336,
(1877)
\bibitem{3}
J. Katriel {\it Bell numbers and coherent states} Physics Letters
A, {\bf 273} (3) (2000): 159-161

\bibitem{4}
M. Schork {\it On the combinatorics of normal ordering bosonic
operators and deformations of it}   J. Phys. A: Math. Gen. {\bf
36} (2003) 4651-4665
\bibitem{5}
J.Katriel, M. Kibler {\it Normal ordering for deformed boson
operators and operator-valued deformed Stirling numbers}
J. Phys. A: Math. Gen.(\bf 25) (1992): 2683-26-91

\bibitem{6}
S.C. Milne {\it A $q$-analog of restricted growth functions,
Dobinski's equality, and Charlier polynomials} ,Trans. Amer. Math.
Soc. {\bf 245} (1978) 89-118

\bibitem{7}
R.Ehrenborg {\it Determinants involving $q$-Stirling numbers},
Advances in Applied Mathematics (\bf 31)(2003): 630-642

\bibitem{8}
Bennett Curtis, Dempsey Kathy J. , Sagan Bruce E.
{\it Partition Lattice $q$-Analogs Related to $q$-Stirling Numbers }
Journal of Algebraic Combinatorics (\bf 03)(3) p.261-283 July 1994

\bibitem{9}
Rajendra S. Deodhar, Murali K. Srinivasan
{\it An inversion number static on set partitions}
preprint  submited to Elsevier Science  12 September 2003

\bibitem{10}
J. Cigler {\it A new $q$-Analogue of Stirling Numbers}
Sitzunber. Abt. II {\bf 201}(1992) : 97-109

\bibitem{11}
A.K.Kwasniewski {\it Poisson, Dobinski, Rota and $q$-coherent states}
ArXiv: math. CO/0402125 v1 9 Feb 2004

\bibitem{12}
A.K.Kwasniewski {\it Main  theorems of extended finite operator
calculus}Integral Transforms and Special Functions, {\bf 14} No 6 (2003): 499-516

\bibitem{13}
A.K.Kwasniewski {\it On Simple Characterizations of Sheffer $psi$-polynomials and
Related Propositions of the Calculus of Sequences} ,Bulletin de la Soc. des Sciences et
de Lettres de Lodz,  {\bf 52}  Ser. Rech. Deform.  36 (2002):45-65, ArXiv: math.CO/0312397 .
\end{thebibliography}
\end{document}